\documentclass[amsmath,amssymb,pra,superscriptaddress,nofootinbib,notitlepage]{revtex4-2}

\newtheorem{definicao}{{\bf Definition}}
\newtheorem{proposicao}{{\bf Proposition}}
\newtheorem{teorema}{{\bf Theorem}}

\usepackage{amssymb}
\usepackage{graphicx}
\usepackage[]{xcolor}

\begin{document}

\title{Impulsive Control on Invariant Surfaces}

\author{C.~C.~Silva Jr.}
\affiliation{Instituto de F\'\i sica, Universidade de Bras\'\i{}lia, CP: 04455, 70919-970 - Bras\'\i{}lia, Brazil}

\author{J.~Mar\~ao}
\affiliation{Universidade Federal do Maranh\~ao, 65085-580 and
Departamento de Matem\'atica e Inform\'atica, Universidade Estadual do Maranh\~ao,
S\~ao Lu\'is - MA, 
650000-001}

\author{A.~Figueiredo}
\affiliation{Instituto de F\'\i sica and International Center for Condensed Matter Physics\\ Universidade de
Bras\'\i{}lia, CP: 04455, 70919-970 - Bras\'\i{}lia, Brazil}

\author{T.~M.~Rocha Filho}
\affiliation{Instituto de F\'\i sica and International Center for Condensed Matter Physics\\ Universidade de
Bras\'\i{}lia, CP: 04455, 70919-970 - Bras\'\i{}lia, Brazil}

\begin{abstract}
An impulsive feedback-adaptive control is developed in order to drive trajectories of a dynamical system
towards an invariant manifold with fixed and spaced impulsive controls.
The approach requires the explicit knowledge of the set of equations defining the invariant
manifold and is based on the concept of stability exponents of invariant manifolds.

\vspace*{1mm}
\end{abstract}

\maketitle

\section{Introduction}

Since the seminal works on synchronization and control of chaotic systems~\cite{winfree,kuramoto,r1,r2} many developments were achieved in the field
(see for instance~\cite{pecorarev,motter} and references therein). The main idea is to constraint
trajectories to converge to an invariant manifold, usually unstable~\cite{alzahrani,artigocontrole,artigoouter},
such as fixed points, periodic orbits, chaotic attractors and hyper-surfaces.
A common approach is to couple an external field to the system to ensure the stability of the
invariant manifold. The external field is said to be continuous if it acts permanently on the system,
or impulsive it acts only for short time intervals or even instantaneously.
The first results in this direction in the literature considered that for the continuous case, the stability condition was equivalent to
the negativity of the Transverse Lyapunov Exponents (TLE) of the invariant manifold~\cite{r3,r4}.
Nevertheless, in many relevant practical situations,
strong desynchronization bursts can occur in the dynamics due to the existence of unstable periodic orbits,
even for the case with negative highest TLE of the chaotic attractor~\cite{r5}. It was also shown from computer simulations
of drive subsystems with a positive TLE, that synchronization of spatiotemporal chaos is possible~\cite{r6}.
Therefore, the negativity of the TLE's is neither a sufficient nor a necessary condition for the stability of a continuous control,
a precise formulation of such conditions being an open question~\cite{artigoleonov}. This difficulty can be overcome by using
feedback adaptive methods~\cite{r7,r8,r9,r10,r11,r12,r13,r14,r15,r16}, with a continuously adapted control,
and stability conditions obtained from Lyapunov Functions.
Other approaches with similar purposes include the switching manifold approach~\cite{r18} 
and the impulsive selective method~\cite{r19}.

For an impulsive control, the stability properties are usually determined
from comparison methods~\cite{r20,r21,r22,r23,r24,r25,r26,r27,r28},
Lyapunov Functions generalized to impulsive systems~\cite{r29,r30,r31,r32,r33,r34} or
an impulsive adaptive-feedback method~\cite{r35}.
Up to our knowledge, the concept of TLE's was not fully exploited for the stability of impulsive controls,
with the additional advantage of allowing unevenly spaced impulsive interventions.
This is the main goal of the present work. Here we use stability exponents that
are closely related to the highest TLE of invariant manifolds~\cite{bainov,samoilenko}.
Our approach can be summarized in two steps: first pull the trajectory back into the direction of the
invariant manifold, and then apply an impulse to correct for instabilities due to the possible existence
of positive TLE's, provided the equations governing the dynamics and the invariant manifold
are explicitly known. Invariant polynomial hyper-surfaces of a polynomial dynamical systems were 
extensively studied in Ref.~\cite{r37}, and a generalization to a class of non-polynomial manifolds of non-polynomial 
vector fields is described in Refs.~\cite{r38,r38b,r40,r40b}.
The method we consider in the present work has a reasonably simple implementation for a control of the convergence speed and stability,
with some advantage over previous approaches~\cite{r39},
with good accuracy in the determination of the intervals between two successive impulses.
As another advantage, it is feedback-adaptive, without requiring prior knowledge on TLE's if the intervals between impulses are upper bounded,
and without loss of control due the presence of unexpected unstable objects inside the invariant manifold.

The paper is organized as follows: in Section~\ref{sec2} we give some basic notations and definitions, in Section~\ref{sec3} we define the concept of Stability Exponent of a given invariant surface and obtain a new way to calculate it, in Section~\ref{sec4} we develop a theory for impulsive control and obtain a general theorem that state the conditions that assure the convergence toward unstable invariant surfaces, in Section~\ref{sec4} we illustrate the results with some examples, in Section~\ref{sec5} we close the article with some concluding remarks.

\section{Basic Notations and Definitions}
\label{sec2}
    
Throughout this work we use the following rules of notation: (a) The real vector space of dimension $m\geq 1$ is denoted by $\mathbb{R}^m$; (b) vectors of a vector space are denoted by bold lower case letters; (c) transformations between two vector spaces are denoted by bold capital letters.
To indicate the components of a vector ${\bf x}\in\mathbb{R}^m$ and  a transformation ${\bf F}:\mathbb{R}^m\rightarrow\mathbb{R}^n$ ($n\leq m$) we use column vectors notation:
\begin{eqnarray}
\label{formula1}
 & & {\bf x}=\left(x_1,\ldots,x_m\right)^t,\;\;x_i\in\mathbb{R},\; i=1,\ldots,m.\nonumber \\
 & & {\bf F}({\bf x})=\left(F_1({\bf x}),\ldots,
	F_n({\bf x})\right)^t,\;\;  F_i:\mathbb{R}^m\rightarrow\mathbb{R},\;  i=1,\ldots,n.
\end{eqnarray}
If all components of a transformation ${\bf F}$ are continuous and differentiable functions of ${\bf x}$, then the transformation is called regular.
Any transformation ${\bf F}:\mathbb{R}^m\rightarrow\mathbb{R}^m$ defines a vector field in $\mathbb{R}^m$.

We can associate to a regular vector field ${\bf F}$ in $\mathbb{R}^m$, the system of $m$ Ordinary Differential Equations (ODE's):
\begin{equation}
\label{formula2}
	\frac{d{\bf x}(t)}{dt}={\bf F}({\bf x}),\; t\in\mathbb{R},\; {\bf x}\in \mathbb{R}^m.
\end{equation}
Solutions of the ODE system in Eq.~(\ref{formula2}) are called trajectories of the vector field ${\bf F}$.

\begin{definicao}(Invariant Set): A set $S\subset\mathbb{R}^m$ is an invariant of a vector field ${\bf F}$ in $\mathbb{R}^m$ if
any of its trajectories that intersect $S$ is contained in $S$. In other words, for each initial condition contained in $S$,
the solution of the ODE system~(\ref{formula2}) is also contained in $S$.
\end{definicao}

\begin{proposicao} Let us consider a regular vector field ${\bf F}$ in $\mathbb{R}^m$ and a regular transformation 
${\bf I}:\mathbb{R}^m\rightarrow\mathbb{R}^p$, where ${\bf I}({\bf x})=(I_1({\bf x}),\ldots,I_p({\bf x}))^t$ with $p\leq m$. If the hyper-surface
$S=\left\{{\bf x}\in \mathbb{R}^m\;|\; {\bf I}({\bf x})=0\right\}$ is an invariant of ${\bf F}$, then ${\bf F}({\bf x})\cdot \nabla_{\bf x}{I}_i({\bf x})=0$  for all ${\bf x}\in S$ and $i=1,\ldots,p$.
\end{proposicao}

The proof of proposition 1 follows directly from to consider the time derivative of $I_i({\bf x}(t))$, where ${\bf x} (t)$ is a trajectory of the vector field ${\bf F}$:
\begin{equation}
\label{formula3}
\displaystyle\left.\frac{d  I_i({\bf x}(t'))}{dt'}\right|_t=\sum_{j=1}^m\left.\frac{\partial I_i{(\bf x)}}{\partial x_j}\right|_t
	\left.\frac{dx_j}{dt'}\right|_t=
{\bf F}({\bf x}(t))\cdot \nabla_{\bf x}{I_i}({\bf x}(t)),\;\; i=1,\ldots,p.
\end{equation}

\section{Semi-Invariants and Stability Exponents}
\label{sec3}

We first state and develop, in more rigorous way, some previous results obtained in Ref.~\cite{marao}. 

\begin{definicao}(Semi-Invariant of dimension $p$):
A transformation ${\bf I}:\mathbb{R}^m\rightarrow\mathbb{R}^p$, with
${\bf I}({\bf x})=(I_1({\bf x}),\ldots,I_p({\bf x}))^t$, is a Semi-Invariant of the vector field ${\bf F}$
if there is a square matrix ${\bf L}({\bf x})=\left[L_{ij}({\bf x})\right]$ of dimension $p$ such that
\begin{equation}
\label{formula3a}
\frac{d {\bf I}(\bf x)}{dt}= {\bf L}({\bf x}){\bf I}({\bf x})\;\equiv\;
\frac{d I_i(\bf x)}{dt}=\sum_{j=1}^{p}{L}_{ij}({\bf x}){I}_j({\bf x} ),\; i=1,\ldots,p,
\end{equation}
where the matrix elements $L_{ij}({\bf x})\in\mathbb{R}$ are regular functions of ${\bf x}$.
\end{definicao}

\begin{definicao}(Non-Singular Semi-Invariant):
A Semi-Invariant ${\bf I}:\mathbb{R}^m\rightarrow\mathbb{R}^p$ ($p\leq m$) of a vector field
${\bf F}$  in $ \mathbb{R}^m$ is non-singular if there is a transformation ${\bf J}:\mathbb{R}^m\rightarrow\mathbb{R}^{m-p}$,
such that the transformation ${\bf T}:\mathbb{R}^m\rightarrow\mathbb{ R}^{m}$, defined as
$${\bf T}({\bf x})=\left(I_1({\bf x}),\ldots,I_p({\bf x}),J_{1}({\bf x}) ,\ldots,J_{m-p}({\bf x})\right)^t$$
is non-singular, i.~e.
$$\det\left(\frac{\partial {\bf T}}{\partial {\bf x}}\right)\neq 0\;\;\forall {\bf x}\in\mathbb{R}^m.$$
\end{definicao}

\begin{proposicao} Lets consider a transformation ${\bf I}:\mathbb{R}^m\rightarrow\mathbb{R}^p$ ($p\leq m$) and
the respective hyper-surface $S=\{{\bf x}\in\mathbb{R}^m\;|\; {\bf I}({\bf x})=0\}$.
If ${\bf I}$ is a non-singular semi-invariant of a vector field ${\bf F}$ in $\mathbb{R}^m$,
then $S$ is a $m-p$ dimensional invariant hyper-surface of ${\bf F}$.
\end{proposicao}

\noindent {\it Proof}: From Definition 2 we have that
\begin{equation}
\label{formula3b}
\forall {\bf x}\in S\;\Rightarrow\; {\bf I}({\bf x})=0\Rightarrow\frac{d{\bf I}({\bf x})}{dt}=0,
\end{equation}
and therefore $S$ is an invariant hyper-surface of ${\bf F}$.
On the other hand, if the semi-invariant ${\bf I}$ is non-singular and using Definition 3 the vectors
$\nabla_{\bf x}I_i({\bf x})$ are linearly independent for all ${\bf x}\in\mathbb{R}^m$,
and define a vector subspace ${\bf E}$ of dimension $p$ contained in the tangent space at ${\bf x}\in\mathbb{R}^m$.
From equation (\ref{formula3}) and (\ref{formula3b}), the tangent vector subspace to the hyper-surface $S$ at ${\bf x}$ is complementary and perpendicular to ${\bf E}$ and has dimension $m-p$. {\it QED}\\

Considering a non-singular semi-invariant as in Definitions 2 and 3, and the associated coordinate transformation ${\bf y}={\bf T}({\bf x})$,
with ${\bf T}$ given in Definition 3, we can rewrite the system of ODE's in (\ref{formula2}) as
$d{\bf y}/dt={\bf G}({\bf y})$, with ${\bf G}:\mathbb{R}^m\rightarrow\mathbb{R}^m$ being a regular vector field such that
\begin{equation}
\label{formula6}
	\frac{d {\bf y}}{dt}={\bf G}({\bf y})\;\equiv\;
\begin{array}{c}        
\displaystyle\frac{d{\bf I}}{dt} = {\bf L}({\bf I},{\bf J}){\bf I},\\ \\
\displaystyle\frac{d{\bf J}}{dt}={\bf P}({\bf I},{\bf J}),
\end{array}
\end{equation}
where ${\bf P}({\bf I},{\bf J})$ is a vector function depending on the choice of the functions ${\bf J}$.

\begin{definicao} (Stability Exponent)
Let us consider a vector field ${\bf F}$ in $\mathbb{R}^m$ with trajectories determined by the solutions of the system of ODEs in
Eq.~(\ref{formula2}), a $p$ dimensional non-singular semi-invariant ${\bf I}$ of ${\bf F}$ and its respective $m-p$ dimensional
invariant hyper-surface $S$.

Let ${\bf x}(t)$ be a trajectory of the ODE system (\ref{formula2}) with initial condition
${\bf x}(0)=\bar{\bf x}_0$. We define
$${\bf I}(t)={\bf I}({\bf x}(t)),\;\;{\bf I}_0={\bf I}(0),\;\;||{\bf I}(t)||=\sqrt{I^2_1(t)+\cdots+I^2_q(t)}.$$
We then compute the following limit for each point ${\bf x}_0$on the invariant surface $S$:
$$
D_S\left({\bf x}_0,{\bf i}_0\right)=\lim_{\bar{\bf x}_0\rightarrow{\bf x}_0}\left[\lim_{t\rightarrow\infty}
\frac{1}{t}\ln\left(\frac{||{\bf I}(t)||}{||{\bf I}_0||}\right)\right],\;\;
{\bf i}_0\equiv\lim_{\bar{\bf x}_0\rightarrow{\bf x}_0} \frac{{\bf I}_0}{||{\bf I}_0||}.
$$ 
If this limit exists, then it is called the stability exponent of the surface $S$ at the point ${\bf x}_0$.
\end{definicao}

The stability exponent for the ODEs system in Eq.~(\ref{formula6}) can be calculated as:
\begin{equation}
\label{formula7}
{\bf J}_0\in \mathbb{R}^{m-p},\;\;{\bf I}_0\in \mathbb{R}^{p},\;\;{\bf i}_0=\lim_{{\bf I}_0\rightarrow 0} \frac{{\bf I}_0}{||{\bf I}_0||}
\;\;\Longrightarrow\;\; D_S\left({\bf J}_0,{\bf i}_0\right)=\lim_{{\bf I}_0\rightarrow 0}\left[\lim_{t\rightarrow\infty}\frac{1}{t}\ln\left(\frac{||{\bf I}(t)||}{||{\bf I}_0||}\right)\right],\;\;
\end{equation}
where
\begin{equation}
\label{formula8}
{\bf I}(t)={\cal I}({\bf I}_0,t),\;\;{\bf J}(t)={\cal J}({\bf J}_0,t),\;\;\;{\bf I}(0)={\bf I}_0, \;\;{\bf J}(0)={\bf J}_0,
\end{equation}
are solutions of Eq.~(\ref{formula6}).

Oseledec theorem~\cite{oseledec,ruelle} establishes the conditions for the existence of the stability exponent in Definition 4.
A method to obtain this exponent was obtained in Ref.~\cite{marao} and is given in Theorem 1 below.
Before to establish this theorem we need the following definitions and proposition.

\begin{definicao}
For a non-singular semi-invariant ${\bf I}$ of a vector field ${\bf F}$ we define the Hermitian matrix:
$${\bf H}({\bf x})=\frac{{\bf L}({\bf x})+{\bf L}^t({\bf x})}{2}, \;\; {\bf x}\in\mathbb{R}^m,$$
where ${\bf L}^t({\bf x})$ is the transpose of ${\bf L}({\bf x})$. The matrix ${\bf H}({\bf x})$ can be written in terms of the variables ${\bf I}$ and ${\bf J}$, which obey the ODEs system in Eq.~\ref{formula6}, as follows:
$$
{\bf H}({\bf I},{\bf J})=\frac{{\bf L}({\bf I},{\bf J})+{\bf L}^t ({\bf I},{\bf J})}{2}.
$$
\end{definicao}

\begin{definicao}Considering the invariant surface $S$ of a non-singular semi-variant ${\bf I}$ as given in Proposition 2,
we denote by ${\bf L}_S({\bf x})$ and ${ \bf H}_S({\bf x})$ the respective matrices ${\bf L}(\bf x)$ and ${\bf H}(\bf x)$
restricted to the domain $ S\subset \mathbb{R}^m$. Considering ${\bf L}(\bf x)$ and ${\bf H}(\bf x)$, written as ${\bf L}({\bf I},{\bf J})$ and 
${\bf H}({\bf I},{\bf J})$ when expressed in terms of the vectors ${\bf I}$ and ${\bf J}$ we have
$${\bf L}_S({\bf J})={\bf L}({\bf 0},{\bf J}),\;\; {\bf H}_S({\bf J})={\bf H}({\bf 0},{\bf J}),\;\; {\bf 0}\in\mathbb{R}^p,\;\; {\bf J}\in\mathbb{R}^{m-p}.$$
\end{definicao}

\begin{proposicao} Let us consider a semi-invariant ${\bf I}$ (non-singular and $p$-dimensional) of a vector field ${\bf F}$ satisfying
Eq.(\ref{formula3a}). The Euclidead norm of ${\bf I}$ and its respective versor ${\bf i}$:
\begin{equation}
\label{formula9}
||{\bf I}||=\sqrt{\sum_{i=1}^{p}I_i^2},\;\; {\bf i}=\frac{{\bf I}}{||{\bf I}||}
\end{equation}
satisfy the equations:
\begin{equation}
\label{formula10}
\frac{d||{\bf I}||}{dt}=||{\bf I}||\langle{\bf i},{\bf H}{\bf i}\rangle,
\;\;\frac{d{\bf i}}{dt}={\bf L}{\bf i}-{\bf i}\langle{\bf i},{\bf H}{\bf i}\rangle,
\end{equation}
where $\langle {\bf A},{\bf B}\rangle$ stands for the canonical dot product between the vectors ${\bf A}$ and ${\bf B}$.
\end{proposicao}
\noindent {\it Proof}: Taking the time derivative of $||{\bf I}||$ with respect to $t$ and using Eq.~(\ref{formula3a}) yields:
\begin{eqnarray}
\frac{d||{\bf I}||}{dt}&=&\frac{1}{2}\frac{1}{||{\bf I}||}\left(\sum_{i=1}^{p}I_i\frac{dI_i}{dt}+\sum_{i=1}^{p}\frac{dI_i}{dt}I_i\right)
=\frac{1}{2}\frac{1}{||{\bf I}||}\left(\sum_{i=1}^{p}I_i\sum_{j=1}^{p}{L}_{ij}{I}_j+\sum_{i=1}^{p}{I}_i\sum_{j=1}^{p}{L}^t_{ij}I_j\right)
\nonumber\\
&=&\frac{1}{||{\bf I}||}\sum_{i=1}^{p}I_i\sum_{j=1}^{p}\frac{\left({L}_{ij}+{L}^t_{ij}\right)}{2}I_j=
\frac{1}{||{\bf I}||}\sum_{i=1}^{p}I_i\sum_{j=1}^{p}{ H}_{ij}I_j
 =  ||{\bf I}||\sum_{i=1}^{p}\frac{I_i}{||{\bf I}||}\sum_{j=1}^{p}{ H}_{ij}\frac{I_j}{||{\bf I}||}
=||{\bf I}||\left<{\bf i},{\bf H}{\bf i}\right>.\nonumber
\end{eqnarray}
In the same way by taking the derivative of ${\bf i}$ we have:
\begin{equation}
\frac{d{\bf i}}{dt}=\frac{d}{dt}\left(\frac{\bf I}{||{\bf I}||}\right)
=\frac{{\bf L}{\bf I}}{||{\bf I}||}-{\bf I}\frac{1}{||{\bf I}||^2}
\frac{d||{\bf I}||}{dt}={\bf L}{\bf i}-\frac{\bf I}{||{\bf I}||^2}||{\bf I}||\left<{\bf i},{\bf H}{\bf i}\right>
={\bf L}{\bf i}-{\bf i}\left<{\bf i},{\bf H}{\bf i}\right>\nonumber
\end{equation}
where we used Definition 2. {\it QED}\\

The latter proposition allows to determine the trajectories of the vector fields ${\bf I}$ and ${\bf J}$ from the system:
\begin{equation}
\label{formula11}
\frac{d||{\bf I}||}{dt}=||{\bf I}||\left<{\bf i},{\bf H}({\bf i}||{ \bf I}||,{\bf J}){\bf i}\right>,\;\;\frac{d{\bf i}}{dt}={\bf L}({\bf i }||{\bf I}||,{\bf J}){\bf i}-{\bf i}\left<{\bf i},{\bf H}({\bf i}|| {\bf I}||,{\bf J}){\bf i}\right>,
\;\;\frac{d {\bf J}}{dt}={\bf P}({\bf i}||{\bf I}||,{\bf J}).
\end{equation}
where we use the relation ${\bf I}={\bf i}||{\bf I}||$. Although this has $m+1$ equations they are not independent.
In fact the system of $p$ equations for ${\bf i}$ has only $p-1$ independent differential equations as $<{\bf i},{\bf i}>=1$
is an invariant scalar.

We now state the first theorem of the present work:

\begin{teorema} A stability exponent of a $m-p$ dimensional invariant surface $S$,
associated with a $p$-dimensional non-singular semi-invariant ${\bf I}$ of a $m$-dimensional vector field
${\bf F }$ ($p\leq m$), can be calculated as
$$ D_S\left({\bf J}_0,{\bf i}_0\right)=\lim_{t\rightarrow\infty}\frac{1}{t}\int_0^{t}\left<{ \bf i}(t'),{\bf H}_S({\bf J}(t')){\bf i}(t')\right>dt',$$
where ${\bf J}(t)$ and ${\bf i}(t)$ constitute a solution to the following initial value problem:
$$\frac{d{\bf i}}{dt}={\bf L}_S({\bf J}){\bf i}-{\bf i}\left<{\bf i},{ \bf H}_S({\bf J}){\bf i}\right>,
\;\;\frac{d {\bf J}}{dt}={\bf P}_S({\bf J }),\;\;{\bf i}(0)={\bf i}_0,\;\;{\bf J}(0)={\bf J}_0.$$
The matrices ${\bf L}_S({\bf J})$ and ${\bf H}_S({\bf J})$ are given in Definition 6. The vector field
${\bf P}_S :\mathbb{R}^{m-p}\rightarrow\mathbb{R}^{m-p}$ is defined by
$${\bf P}_S({\bf J})\equiv{\bf P}(0,{\bf J}),\;\;0\in\mathbb{R}^p,\;\; {\bf J}\in\mathbb{R}^{m-p},$$
for ${\bf P}({\bf I},{\bf J})$ given in Eq.~(\ref{formula6}).
\end{teorema}

\noindent{\it Proof:} From Eqs.~(\ref{formula7}), (\ref{formula8}) and~(\ref{formula11}) we have:
$$ D_S\left({\bf J}_0,{\bf i}_0\right)=\lim_{{\bf I}_0\rightarrow 0}\left[\lim_{t\rightarrow\infty}\frac {1}{t}\ln\left(\frac{||{\bf I}(t)||}{||{\bf I}_0||}\right)\right]=
\lim_{{\bf I}_0\rightarrow 0}\left[\lim_{t\rightarrow\infty}\frac{1}{t}\int_0^{t}
\langle{\bf i},{\bf H}({\bf i}||{\bf I}||,{\bf J}){\bf i}\rangle_{t'}\,dt^\prime\right],$$
where the subscript $t^\prime$ means that all functions in the dot product are evaluated for $t^\prime\in[0,t]$.
As $S$ is an invariant surface we also have that
\begin{equation}
\label{formula12}
\lim_{{\bf I}_0\rightarrow 0}{\bf I}(t')=\lim_{{\bf I}_0\rightarrow 0}{\cal I}({\bf I}_0,t ')=0\;\;\Rightarrow\;\;\lim_{{\bf I}_0\rightarrow 0}{||\bf I}(t')||=0,\;\;\forall t '\in[0,t].
\end{equation}
Therefore,  $D_S$ can be written as:
\begin{equation}
\label{formula13}
D_S\left({\bf J}_0,{\bf i}_0\right)=\lim_{||{\bf I}(t')||\rightarrow 0}\left[\lim_{t\rightarrow \infty}\frac{1}{t}\int_0^{t}\left<{\bf i},{\bf H}({\bf i}||{\bf I}||,{\bf J}){\bf i}\right>_{t'}\,dt'\right]=
\lim_{t\rightarrow\infty}\frac{1}{t}\int_0^{t}\left<{\bf i},{\bf H}(0,{\bf J}){\bf i }\right>_{t'}\,dt'.
\end{equation}
Plugging Eq.~(\ref{formula12}) into Eq.~(\ref{formula11}) we obtain the system of ODEs:
\begin{eqnarray}
\label{formula14}
&&\lim_{||{\bf I}||\rightarrow 0}\left[\frac{d{\bf i}}{dt}={\bf L}({\bf i}||{\bf I}||,{\bf J}){\bf i}-{\bf i}\left<{\bf i},{\bf H}({\bf i}||{\bf I}| |,{\bf J}){\bf i}\right>\right],
\;\;\lim_{||{\bf I}||\rightarrow 0}\left[\frac{d {\bf J}}{dt}={\bf P}({\bf i}|| {\bf I}||,{\bf J})\right]\nonumber \\
&&\Rightarrow\;\; \frac{d{\bf i}}{dt}={\bf L}(0,{\bf J}){\bf i}-{\bf i}\left<{\bf i},{\bf H}(0,{\bf J}){\bf i}\right>,\;\;\frac{d {\bf J}}{dt}={\bf P}(0,{\bf J}),
\end{eqnarray}
with initial condition:
\begin{equation}
\label{formula15}
{\bf i}(0)={\bf i}_0=\lim_{{\bf I}_0\rightarrow 0}\frac{{\bf I}_0}{||{\bf I}_0|| },\;\;{\bf J}(0)={\bf J}_0.
\end{equation}
Taking into account the definitions of ${\bf L}_S({\bf J})$, ${\bf H}_S({\bf J})$ and ${\bf P}_S({\bf J})$, then equations (\ref{formula13}), (\ref{formula14}) and (\ref{formula15}) represent the result stated in the theorem. {\it QED}\\

The stability exponent can therefore be determined as a time average along trajectories in the invariant hyper-surface $S$.
The value of the exponent is then determined from the attractor objects in $S$ and not by the initial condition ${\bf J}(0)={\bf J}_0$.
Indeed, Oseledet theorem establishes that a stability exponent associated with a given attractor
must depend only on the initial value ${\bf i}(0)={\bf i }_0$. 

\section{Impulsive Systems and Convergence Control}
\label{sec4}

In this section we define an impulsive system coupled with a vector field
and present some results to control the convergence of their trajectories  to an unstable invariant surfaces, associated to
a semi-invariant of the vector field.  In what follows we consider:
\begin{enumerate}
\item A $m$-dimensional vector field ${\bf F}$ with trajectories that are solutions of  Eq.~(\ref{formula2}).
\item A $p$-dimensional non-singular semi-invariant ${\bf I}$ of ${\bf F}$ with $p\leq m$ satisfying Eq.~(\ref{formula3a}).
\item The $m-p$ dimensional invariant surface $S$ associated to the semi-invariant ${\bf I}$.
\item The transformation of coordinates defined by the vector field ${\bf T}$ and the respective vector field ${\bf G}$
as given in Definition 3, expressed in terms of the vectors ${\bf I}$ and ${\bf J}$,
with trajectories satisfying Eq.~(\ref{formula6}).
\end{enumerate}

\begin{definicao}
(Impulsive System with Initial Condition):
An impulsive system with an initial condition at a certain time $t_0\in R$, associated with the vector field ${\bf F}$, is defined by the following set of equations:
\begin{eqnarray}
\label{formula16}
\begin{array}{l}
\displaystyle \frac{d{\bf x}(t)}{dt}={\bf F}({\bf x}(t)), \;\;\;\;\forall t\neq t_n, \\ \\
\Delta{\bf x}(t)={\bf x}(t_n^{+})-{\bf x}(t_n^{-})={\bf G}_n,\;\;\;\;\forall t=t_n,\\ \\
{\bf x}(t_0)={\bf x}_0\in\mathbb{R}^m,\;\;\;\;t_0<t_n,
\end{array}
\end{eqnarray}
where $t_n^+$ and $t_n^-$ respectively denote the left and right limits for $t\rightarrow t_n$. The sequence
$t_n$ for $n\in\mathbb{N}^+$ is such that $t_{n-1}<t_{n}$ and $\lim_{n\rightarrow\infty}t_n=\infty$. $\Delta {\bf x}(t)$ is called impulse and
$\Delta_n=t_{n}-t_{n-1}$ denotes the time interval between two impulses.
The impulsive vector field is given by a sequence of vectors ${\bf G}_n=(G_{n1},\ldots,G_{nm})^t\in\mathbb{R}^m$ with
$n\in\mathbb{N}^+$.

In terms of the vectors ${\bf I}$ and ${\bf J}$ the impulsive system is written as
\begin{eqnarray}
\label{formula17}
\begin{array}{l}
\displaystyle \frac{d {\bf I}(t)}{dt}= {\bf L}({\bf I}(t),{\bf J}(t)){\bf I}(t),
\;\; \frac{d {\bf J}(t)}{dt}={\bf P}({\bf I}(t),{\bf J}(t)), \;\;\;\;\forall t\neq t_n, \\ \\
\Delta{\bf I}(t)={\bf I}(t_n^{+})-{\bf I}(t_n^{-})={\bf G}^I_n,\;\;\Delta{\bf J}(t)={\bf J}(t_n^{+})-{\bf J}(t_n^{-})={\bf G}^J_n, \;\;\;\;\forall t=t_n, \\ \\
{\bf I}(t_0)={\bf I}_0\in\mathbb{R}^p,\;\;{\bf J}(t_0)={\bf J}_0\in\mathbb{R}^{m-p},\;\;\;\;t_0<t_n,
\end{array}
\end{eqnarray}
where for every $n\in\mathbb{N}$, we have ${\bf G}^I_n=(G^I_{n1},\ldots,G^I_{np})^t\in\mathbb{ R}^ p$
and ${\bf G}^J_n=(G^I_{n1},\ldots,G^I_{nm-p})^t\in\mathbb{R}^{m-p}$.
The sequence of vectors ${\bf G}^I_n$ and ${\bf G}^J_n$ are unequivocally determined by the sequence of vectors ${\bf G}_n$
and the coordinate transformation defined by ${\bf T}$.
\end{definicao}
\begin{definicao}
Considering a solution ${\bf x}(t)$ of the impulsive system in (\ref{formula16}), or alternatively,
a solution $[{\bf I}(t),{\bf J}(t)]$ of the system (\ref{formula17}), we define the exponents $\beta_n[t]$ :
\begin{equation}
\label{formula18}
\begin{array}{l}
\displaystyle||{\bf I}(t)||=||{\bf I}(t_{n-1}^{+})||\exp(\beta_n[t])\;\;\Leftrightarrow\;\; \beta_n[t]\equiv\ln\left(\frac{||{\bf I}(t)||}{||{\bf I}(t_{n-1}^{+})||}\right),\;\;n\in\mathbb{N},\;\; t_{n-1}<t<t_n,\\
\displaystyle||{\bf I}(t_{n}^-)||=||{\bf I}(t_{n-1}^{+})||\exp(\beta_n[t_n])\;\;\Leftrightarrow\;\; \beta_n[t_n]\equiv\ln\left(\frac{||{\bf I}(t_n^-)||}{||{\bf I}(t_{n-1}^{+})||}\right),\;\;n\in\mathbb{N},
\end{array}
\end{equation}
with $t_0^+=t_0$. From Eq.~(\ref{formula11}) we have:
\begin{equation}
\label{formula19}
\begin{array}{l}
\displaystyle\beta_n[t]=\int_{t_{n-1}^+}^t \left<{\bf i},{\bf H}({\bf i}||{\bf I}||,{\bf J}){\bf i}\right>dt',\;\;n\in\mathbb{N},\;\; t_{n-1}<t<t_n,\\
\displaystyle\beta_n[t_n]=\int_{t_{n-1}^+}^{t_n^-} \left<{\bf i},{\bf H}({\bf i}||{\bf I}||,{\bf J}){\bf i}\right>dt',\;\;n\in\mathbb{N}.
\end{array}
\end{equation}
\end{definicao}
\begin{definicao}Considering a solution ${\bf x}(t)$ of the impulsive system in (\ref{formula16}), or alternatively, the solution $[{\bf I}(t),{\bf J}(t)]$ of the system (\ref{formula17}), we define the exponents $A_n$:
\begin{equation}
\label{formula20}
||{\bf I}(t_n^{+})||=||{\bf I}(t_n^{-})+{\bf G}^I_n||=||{\bf I}(t_n^{-})||\exp(A_n)
\;\;\Leftrightarrow\;\; A_n\equiv\ln\left(\frac{||{\bf I}(t_n^{+})||}{||{\bf I}(t_n^{-})||}\right),\;\;n\in\mathbb{N}.
\end{equation}
\end{definicao}
\begin{definicao} From the exponents $A_n$ and $\beta_n$ in Definitions 8 and 9, we can define the exponents $B_n$:
\begin{equation}
\label{formula21}
B_n=A_n+\beta_n[t_n]\;\;\Leftrightarrow\;\; A_n=B_n-\beta_n[t_n].
\end{equation}
\end{definicao}

Using the exponents in Definitions~8, 9 and~10 we determine the limit of $\lim_{t\rightarrow\infty}||{\bf I}(t)||$ for a solution ${ \bf x}(t)$,
respectively $[{\bf I}(t),{\bf J}(t)]$, of the impulsive system in Definition 7.
From Definitions~8 and~9 we obtain the recurrence relation:
\begin{equation}
\label{formula22}
||{\bf I}(t_{n}^{+})||=||{\bf I}(t_{n-1}^{+})||\exp(A_{n}+\beta_{n}[t_{n}]).
\end{equation}
Iterating the recurrence relation (\ref{formula22}) up to a certain value $n>n_0$ we obtain:
\begin{equation}
\label{formula23}
||{\bf I}(t_{n}^{+})||=||{\bf I}(t_{n_0-1}^{+})||\exp\left(\sum_{p=n_0}^{n}A_q+\sum_{p=n_0}^{n}\beta_q[t_q]\right), \;\; n>n_0.
\end{equation}
Considering $\beta_n[t]$ defined in (\ref{formula18}) and $||{\bf I}(t_{n}^{+})||$ in (\ref{formula23}), it can be show that
\begin{equation}
\label{formula24}
||{\bf I}(t)||=||{\bf I}(t_{n_0-1}^{+})||\exp\left(\sum_{p=n_0}^{n}A_q+\sum_{p=n_0}^{n}\beta_q[t_q]+\beta_{n+1}[t]\right), \;\; n>n_0,\;\; t_n\leq t <t_{n+1}.
\end{equation}
Integrating the differential equation for $||{\bf I}||$ in (\ref{formula11}) we obtain:
\begin{eqnarray}
\label{formula25}
\hspace*{-5mm}
\sum_{p=n_0}^{n}\beta_q[t_q]+\beta_{n+1}[t]&=&\sum_{p=n_0}^{n}\int_{t_{p-1}^+}^{t_q^-}
\langle{\bf i},{\bf H}({\bf i}||{\bf I}||,{\bf J}){\bf i}\rangle dt'
+\int_{t_{n}^+}^t \langle{\bf i},{\bf H}({\bf i}||{\bf I}||,{\bf J}){\bf i}\rangle dt'
=\int_{t_{n_0-1}}^{t} \langle{\bf i},{\bf H}({\bf i}||{\bf I}||,{\bf J}){\bf i}\rangle dt',
\end{eqnarray}
for $t_n\leq t <t_{n+1}$.
Similarly, from the definition $A_n$ in Eq.~(\ref{formula21}), we also have:
\begin{eqnarray}
\label{formula26}
\hspace*{-5mm}
\sum_{p=n_0}^{n}A_q&=&\sum_{p=n_0}^{n}B_q-\sum_{p=n_0}^{p}\beta_q[t_q]=\sum_{p=n_0}^{n}B_q-\sum_{p=n_0}^{n}\int_{t_{p-1}^+}^{t_q^-}
\langle{\bf i},{\bf H}({\bf i}||{\bf I}||,{\bf J}){\bf i}\rangle dt
=\sum_{p=n_0}^{n}B_q-\int_{t_{n_0-1}}^{t_n} \langle{\bf i},{\bf H}({\bf i}||{\bf I}||,{\bf J}){\bf i}\rangle dt'.
\end{eqnarray}
Assuming that the following limit exists:
\begin{equation}
\label{formula27}
\omega=\lim_{t\rightarrow\infty}\frac{\displaystyle\int_{t_{n_0-1}}^{t}
\langle{\bf i},{\bf H}({\bf i}||{\bf I}||,{\bf J}){\bf i}\rangle dt'}{(t-t_{n_0-1})}
=\lim_{t_n\rightarrow\infty}\frac{\displaystyle\int_{t_{n_0-1}}^{t_n}
\langle{\bf i},{\bf H}({\bf i}||{\bf I}||,{\bf J}){\bf i}\rangle dt'}{(t_n-t_{n_0-1})},
\end{equation}
then we obtain from Eqs.~(\ref{formula25}), (\ref{formula26}) and~(\ref{formula27}):
\begin{equation}
\label{formula28}
\lim_{n,t\rightarrow\infty}\left(\sum_{p=n_0}^{n}\beta_q[t_q]+\beta_{n+1}[t]\right)=\lim_{t\rightarrow\infty}[(t-t_{n_0-1})\omega],
\;\;\lim_{n\rightarrow\infty}\sum_{p=n_0}^{n}A_q=\sum_{p=n_0}^{\infty}B_q-\lim_{t_n\rightarrow\infty}[(t_n-t_{n_0-1})\omega].
\end{equation}
From Eqs.~(\ref{formula23}) and~(\ref{formula28}) we have:
\begin{eqnarray}
\label{formula29}
\lim_{t\rightarrow\infty}||{\bf I}(t)||&=&||{\bf I}(t_{n_0-1}^{+})||\lim_{t\rightarrow\infty}\exp\left(\sum_{p=n_0}^{n}A_q+\sum_{p=n_0}^{n}\beta_q[t_q]+\beta_{n+1}[t]\right)\nonumber\\
&=&||{\bf I}(t_{n_0-1}^{+})||\exp\left(\sum_{p=n_0}^{\infty}B_q-\lim_{t_n\rightarrow\infty}[(t_n-t_{n_0-1})\omega]+\lim_{t\rightarrow\infty}[(t-t_{n_0-1})\omega]\right)\nonumber\\
&=&||{\bf I}(t_{n_0-1}^{+})||\exp\left(\sum_{p=n_0}^{\infty}B_q-\lim_{t_n\rightarrow\infty}[(t-t_{n})\omega]\right),\;\; t_n\leq t <t_{n+1}
\end{eqnarray}
We then obtain, using the recurrence relation in Eq.~(\ref{formula22}), that
\begin{equation}
\label{formula30}
||{\bf I}(t_{n_0-1}^{+})||=||{\bf I}_0||\exp\left(\sum_{p=n_0}^{n_0-1}B_q\right),\;\;{\bf I}_0={\bf I}(t_0)={\bf I}({\bf x}(t_0)),\;\;t_0=t_0^+.
\end{equation}
Plugging Eq.~(\ref{formula30}) into Eq.~(\ref{formula29}) yields
\begin{equation}
\label{formula31}
\lim_{t\rightarrow\infty}||{\bf I}(t)||=||{\bf I}_0||\exp\left(\sum_{q=1}^{\infty}B_q-\lim_{t_n\rightarrow\infty}[(t-t_{n})\omega]\right),\;\; t_n\leq t <t_{n+1}.
\end{equation}

We are now in a position to state the main result of the present work, establishing conditions for the convergence of
the impulsive system solution to an invariant hyper-surface $S$ associated to a semi-invariant ${\bf I}$ of a vector field ${ \bf F}$:

\begin{teorema}
Considering a solution ${\bf x}(t)$ of the impulsive system in (\ref{formula16}), or alternatively, the solution $[{\bf I}(t),{\bf J}(t)]$ of the system (\ref{formula17}), then
\begin{equation}
\label{formula32}
\lim_{t\rightarrow\infty}{\bf x}(t)\in S\;\left(\lim_{t\rightarrow\infty}||{\bf I}(t)||=0\right)\;\;\Leftrightarrow\;\;\sum_{q=1}^{\infty}B_q-\lim_{t_n\rightarrow\infty}[(t-t_{n})D_S]=-\infty,\;\; t_n\leq t <t_{n+1},
\end{equation}
where $D_S$ is a stability exponent of the surface $S$, according to definition 3 (respectively calculated in theorem 1).
\end{teorema}
\noindent{\it Proof}: From equation (\ref{formula31}) we obtain
\begin{equation}
\label{formula33}
\lim_{t\rightarrow\infty}||{\bf I}(t)||=0\;\;\Leftrightarrow\;\;\sum_{q=1}^{\infty}B_q-\lim_{ t_n\rightarrow\infty}[(t-t_{n})\omega]=-\infty,\;\; t_n\leq t <t_{n+1}.
\end{equation}
In particular, for $t=t_n$ we have
$$
\lim_{t\rightarrow\infty}||{\bf I}(t_n^+)||=0\;\;\Leftrightarrow\;\;\sum_{q=1}^{\infty}B_q=- \infty.
$$
From the definition of $\omega$ in Eq.~(\ref{formula27}) and the expression for $D_S$ in Theorem 1, one can show that
$$
\lim ||{\bf I}(t_{n_0-1}^+)||=0\;\;\Rightarrow\;\;
\omega=\lim_{t\rightarrow\infty}\frac{\displaystyle\int_{t_{n_0-1}}^{t} \left<{\bf i},{\bf H}(0,{\bf J}){\bf i}\right>dt'}{(t-t_{n_0-1})}=\frac{\displaystyle\int_{t_{n_0-1}}^{t} {\langle\bf i},{\bf H}_S({\bf J}){\bf i}\rangle dt'}{(t-t_{n_0-1})}=D_S
$$
Thus, taking the limit $n_0\rightarrow\infty$ in Eqs.~(\ref{formula29}) and~(\ref{formula30}) is equivalent to
considering the limit $\lim_{n_0\rightarrow\infty}\omega=D_S$ in Eq.~(\ref{formula31}),
and thence Eq.~(\ref{formula33}) implies Eq.~(\ref{formula32}). {\it QED}\\

Theorem 2 allows us to formulate some conditions on the exponent $B_n$ and on the intervals $\Delta_n$ to assure convergence
of the impulsive system trajectories to the invariant surface $S$.

\begin{proposicao}
If there is a real number $\delta>0$ such that $\Delta_n=t_n-t_{n-1}\leq\delta$ for all $n\in\mathbb{N}$, then
\begin{equation}
\label{formula34}
\sum_{q=1}^{\infty}B_q=-\infty\;\;\Rightarrow\;\;\lim_{t\rightarrow\infty}{\bf x}(t)\in S\;\left(\lim_{t\rightarrow\infty}||{\bf I}(t)||=0\right).
\end{equation}
\end{proposicao}

\begin{proposicao}
If there is a positive real number $M>0$ such that $D_S\leq M$, then
\begin{equation}
\label{formula35}
\sum_{q=1}^{\infty}B_q+M\lim_{n\rightarrow\infty}\Delta_{n+1}=-\infty\;\;\Rightarrow\;\;\lim_{t\rightarrow\infty}{\bf x}(t)\in S\;\left(\lim_{t\rightarrow\infty}||{\bf I}(t)||=0\right).
\end{equation}
\end{proposicao}

\begin{proposicao}
Consider that there is a real number $M>0$ such that $D_S\leq M$. Suppose the exponent in Definition 10
is given by $B_n=-\alpha\Delta_n$ with $\alpha>0$, then
\begin{equation}
\label{formula36}
\lim_{n\rightarrow\infty}\left(-\alpha(t_n-t_0)+M\Delta_{n+1}\right)=-\infty\;\;\Rightarrow\;\;\lim_{t \rightarrow\infty}{\bf x}(t)\in S\;\left(\lim_{t\rightarrow\infty}||{\bf I}(t)||=0\right).
\end{equation}
In particular, for any real number $C$ and  $\kappa>0$ such that $0<\kappa<\alpha$, we have
\begin{equation}
\label{formula37}
\Delta_{n+1}\leq\frac{\beta(t_n-t_0)}{M}+C,\;\;\forall n\geq n',\;\;n'\geq 1\;\;\Rightarrow\;\;\lim_{t\rightarrow\infty}{\bf x}(t)\in S\;\left(\lim_{t\rightarrow\infty}||{\bf I}(t)||=0\right).
\end{equation}
\end{proposicao}
Propositions~4 and~5 follow straightforwardly from Theorem~2, and Proposition~6 follows from Proposition 5.

It is important to highlight that Theorem~2, as well as Propositions~4, 5 and~6, assume the existence of the stability exponent
for the asymptotic behavior of the trajectories contained in the invariant hyper-surface $S$.
However, Proposition~5 can be reformulated without this hypothesis by considering that the trajectories of the impulsive system
in Definition~7 are contained in a compact set of $\mathbb{R}^m$.
 
Equation~(\ref{formula24}) for $n_0=1$ leads to
\begin{equation}
\label{formula38}
||{\bf I}(t)||=||{\bf I}_0||\exp\left(\sum_{q=1}^{n}A_q+\sum_{p=n_0}^{n}\beta_q[t_q]+\beta_{n+1}[t]\right), \;\; n\geq 1,\;\; t_n\leq t <t_{n+1}.
\end{equation}
Now using the values of $A_n$ and $\beta_{n+1}[t]$ from Eqs.~(\ref{formula21}) and~(\ref{formula19}) we obtain
\begin{equation}
\label{formula39}
||{\bf I}(t)||=||{\bf I}_0||\exp\left(\sum_{q=1}^{n}B_q+\int_{t_{n}^+}^t \left<{\bf i}({\bf x}),{\bf H}({\bf x}){\bf i}({\bf x})\right>dt'\right), \;\; n\geq 1,\;\; t_n\leq t <t_{n+1},
\end{equation} 
for any solution ${\bf x}(t)$ of the impulsive system in Definition 7.

As the matrix ${\bf H}({\bf x})$ is Hermitian, there is an orthonormal basis formed by the eigenvalues of ${\bf H}({\bf x})$:
${\bf v}_1({\bf x}),\ldots,{\bf v_p}(\bf x)$, with respective eigenvalues $\lambda_1({\bf x}),\ldots,\lambda_p({\bf x})$,
such that ${\bf i}({\bf x})=c_1({\bf x}){\bf v}_1({\bf x})+ \ldots+c_p({\bf x}){\bf v}_p({\bf x})$, with $c_1^2({\bf x})+\ldots+c_p^2({\bf x })=1$.
Therefore, the integral term in Eq.~(\ref{formula39}) can written as
$$\int_{t_{n}^+}^t \left<{\bf i}(t'),{\bf H}(t'){\bf i}(t')\right>dt' = \int_{t_n^+}^{t}\left(c_1^2(t')\lambda_1(t')
+\ldots+c_p^2(t')\lambda_p(t')\right)dt',$$
with ${\bf H}(t)\equiv{\bf H}({\bf x}(t))$, $c_i(t)\equiv c_i({\bf x}(t))$ and $\lambda_i(t)\equiv\lambda_i({\bf x}(t))$,
for a given solution ${\bf x}(t)$ of the impulsive system.

Denoting by $\lambda_H({\bf x})$ the largest eigenvalue of ${\bf H}({\bf x})$ for ${\bf x}\in\mathbb{R}^m$, one can show that
$$\int_{t_{n}^+}^t \left<{\bf i}(t'),{\bf H}(t'){\bf i}(t')\right>dt' \leq\int_{t_n^+}^t\lambda_H({\bf x}(t'))dt',$$
which from Eq.~(\ref{formula39}) implies the following inequality:
\begin{equation}
\label{formula40}
||{\bf I}(t)||\leq||{\bf I}_0||\exp\left(\sum_{q=1}^{n}B_q+\int_{t_n^+}^t \lambda_H({\bf x}(t'))dt\right), \;\; n\geq 1,\;\; t_n\leq t <t_{n+1}.
\end{equation}
We then state our next theorem:

\begin{teorema}
Consider a trajectory ${\bf x}(t)$ of the impulsive system given in Definition~7.
Let us suppose the existence of a compact set $U\subset\mathbb{R}^m$ such that ${\bf x}(t) \subset U$ for all
$t\geq t_0$. Then $\displaystyle M=\max_{{\bf x}\in U}\lambda_H({\bf x})$ is a finite number and the following condition holds:
\begin{equation}
\label{formula41}
\sum_{q=1}^{\infty}B_q+M\lim_{n\rightarrow\infty}\Delta_{n+1}=-\infty\;\;\Rightarrow
\;\;\lim_{t\rightarrow\infty}{\bf x}(t)\in S\;\left(\lim_{t\rightarrow\infty}||{\bf I}(t)||=0\right).
\end{equation}
\end{teorema}
\noindent{\it Proof:} From Eq.~(\ref{formula40}) and the definition of $M$ above, the following inequality holds:
\begin{eqnarray}
||{\bf I}(t)||&\leq&||{\bf I}_0||\exp\left(\sum_{q=1}^{n}B_q+\int_{t_n^+}^ t\lambda_H({\bf x}(t'))dt\right)
\;\leq\; ||{\bf I}_0||\exp\left(\sum_{q=1}^{n}B_q+\int_{t_n^+}^t M\,dt'\right)\nonumber\\
&\leq& ||{\bf I}_0||\exp\left(\sum_{q=1}^{n}B_q+ M(t-t_n)\right)\;\leq\; ||{\bf I}_0||\exp\left(\sum_{q=1}^{n}B_q+ M(t_{n+1}-t_n)\right), \;\; n\geq 1,\;\; t_n\leq t <t_{n+1}. \nonumber
\end{eqnarray}
Remembering that $\Delta_{n+1}=t_{n+1}-t_n$ and taking the limit $n\rightarrow\infty$, we get
$$\lim_{n\rightarrow\infty}||{\bf I}(t)||\leq||{\bf I}_0||\exp\left(\sum_{q=1}^{\infty}B_q+ M\lim_{n\rightarrow\infty}\Delta_{n+1}\right),$$
from which the relation stated in Eq.~(\ref{formula41}) follows.
We note that $M$ is a finite number due to the fact that the function $\lambda_H({\bf x})$ is a regular function on a compact set. {\it QED}\\

\section{Applications}
\label{sec5}

\subsection{The Lorenz system: Convergence towards a Fixed Point}

As a first illustration of our approach, we consider the Lorenz system~\cite{lorenz}:
\begin{eqnarray}
\label{exA1}
 & & \frac{dx_{1}}{dt}=\sigma(x_{2}-x_{1}),\nonumber\\
 & & \frac{dx_{2}}{dt}=rx_{1}-x_{2}-x_{1}x_{3},\nonumber\\
 & & \frac{dx_{3}}{dt}=-bx_{3}+x_{1}x_{2}.
\end{eqnarray}
The fixed point at the origin $x_1=x_2=x_3=0$ is unstable for $r>1$.
We now use our approach to show how to drive solutions of the system towards this fixed point by a proper choice of impulsive interventions.
For that purpose, we consider as invariant manifold:
\begin{equation}
S=\{{\bf x}\in\mathbb{R}^3\;|\;x_1=0,\;x_2=0,\;x_3=0\},
	\label{defSap1}
\end{equation}
with ${\bf I}=(x_1,x_2,x_3)^T$. The matrices ${\bf L}$ and ${\bf L}_S$,
as defined in Eqs.~(\ref{formula3a}) and Definition~5, respectively, are given by:
$$
{\bf L}(t)=  \left[
\begin{array}{ccc}
 - \sigma  & \sigma  & 0 \\
r & -1 &  - x_1(t) \\
0 & x_1(t) &  - b
\end{array}
\right],\;\;\;\;
{{\bf L}_{S}}= \left[
{\begin{array}{ccc}
 - \sigma  & \sigma  & 0 \\
r & -1 & 0 \\
0 & 0 &  - b
\end{array}}
 \right].\;
$$
To define the impulsive vector field, we follow the prescriptions of Sec.~\ref{sec4} and chose ${\bf G}_n$ in Eq.~(\ref{formula16}) and $A_n$ in Eq.~(\ref{formula20}) such that
\begin{equation}
{\bf x}(t_n^+)=\exp(A_n){\bf x}(t_n^-),\;\;\;A_n=-\ln\frac{||x(t_n^-)||}{||x(t_{n-1}^+)||}-\alpha\Delta_n\;\;(\Delta_n=t_n-t_{n-1}).
\label{esquecida1}
\end{equation}
Then, according to Eqs.~(\ref{formula18}) and (\ref{formula21}) we have that $B_n=-\alpha\Delta_n$. The sequence of times $t_n$, when the impulses occur, is defined as
\begin{equation}
\label{esquecida2}
t_{n+1}=t_n+\frac{\kappa}{D_S}(t_n-t_0),\;\;t_1>t_0\;\;(n=1,2,\ldots,\infty),
\end{equation}
where $t_0$ is the initial time for a given initial condition ${\bf x}(t_0)={\bf x}_0$ and $\Delta_1=t_1-t_0$, the time interval for the first impulse to occur, is arbitrary.

For the present example, we take the parameter values $b=8/3$, $\sigma=10$, $r=28$
such that the system is chaotic~\cite{r41}, and chose $\alpha=5$ in Eq.~(\ref{esquecida1}),
$t_0=0$, $t_1=\Delta_1=0.01$.
The stability exponent is obtained using Theorem 1 and is given by the largest eigenvalue of $L_S$ and given by
$D_S=-11/2+\sqrt{1201}/2>0$.

\begin{figure}[ht]
	\includegraphics[scale=0.3]{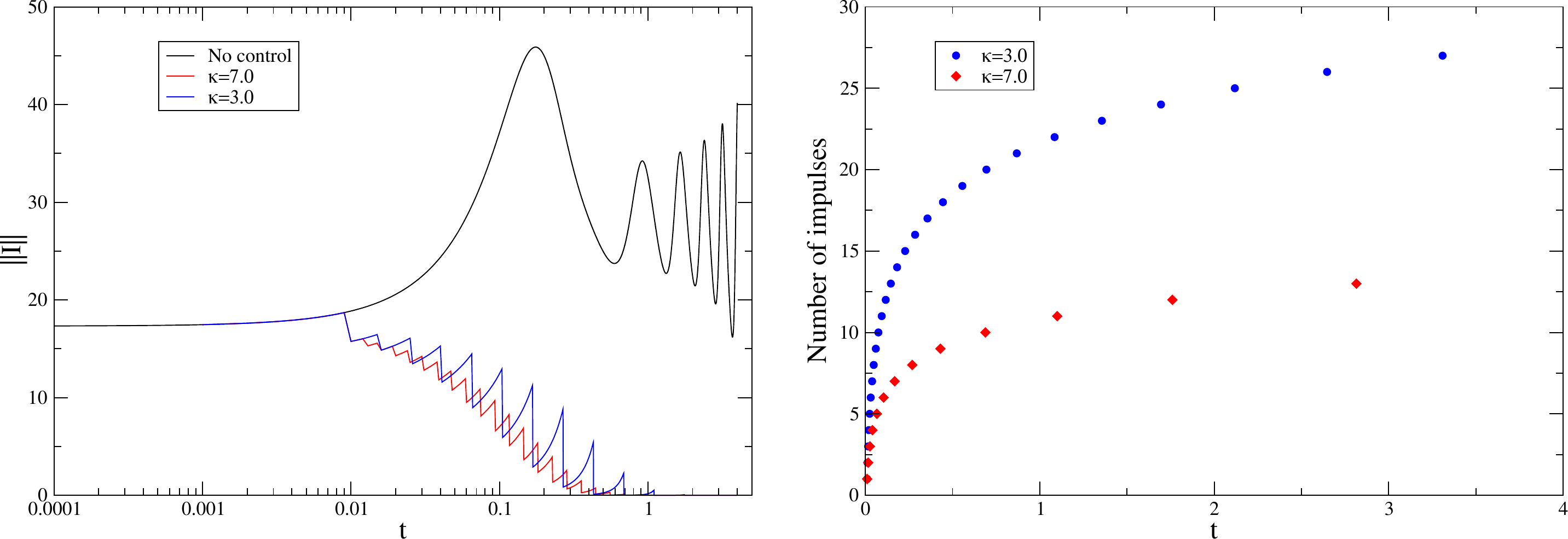}
\caption{Left panel: time evolution $||{\bf I}(t)||$ for the Lorenz system with no control (without impulses)
and two different impulsive controls with $\kappa=3.0$ and $\kappa=7.0$.
Right panel: cumulative number of impulses as a function of time $t$ for $\kappa=3.0$ and $\kappa=7.0$.}
\label{fig1}
\end{figure}

The left panel of Fig.~\ref{fig1} shows the time evolution of $||{\bf I}(t)||$ from a numeric integration of the Lorenz system using a
fourth order Runge-Kutta integrator, with time step $\Delta t=0.001$. We considered the cases with $\beta=3.0$ and $\beta=7.0$.
For comparison purposes we show also the case with no control.
For both impulsive controls the solution rapidly converges to the fixed point.
The cumulative number of impulses is shown in the right panel of Fig.~\ref{fig1}, with frequency of impulses decreasing with time.

Let us observe that from Eq.~(\ref{esquecida2}) we have that $\Delta_n=t_{n+1}-t_n=\kappa/D_S(t_n-t_0)$, which shows that proposition 6 establishes a sufficient (case $\kappa=3<\alpha=5$ in Fig.~\ref{fig1}) but not a necessary (case $\kappa=7>\alpha=5$ in Fig.~\ref{fig1}) condition for the convergence of the impulsive control defined by Eqs.~(\ref{esquecida1}) and (\ref{esquecida2}).

\subsection{Synchronization of two Lorenz systems}

Let us consider the dynamical system given by the Lorenz system in Eq.~(\ref{exA1}) 
and the following Lorenz system with an additional linear coupling term:
\begin{eqnarray}
\label{C1}
 & & \frac{dy_{1}}{dt}=\sigma(y_{2}-y_{1})+c(x_{1}-y_{1}),\nonumber\\
 & & \frac{dy_{2}}{dt}=ry_{1}-y_{2}-y_{1}y_{3},\nonumber\\
 & & \frac{dy_{3}}{dt}=-by_{3}+y_{1}y_{2},
\end{eqnarray}
where the values of the parameter $\sigma$, $r$ and $b$ are the same as in Eq~(\ref{exA1}) and $c$ is the coupling constant. 

Let us define ${\bf x}\equiv(x_{1},x_{2},x_{3})^t$ and ${\bf y}\equiv(y_{1},y_{2},y_{3})^t$.
In order to synchronize the dynamics of $\bf x$ and $\bf y$, we consider the three-dimensional hyper-surface $S$ defined by
$I_1=x_1-y_1=0$, $I_2=x_2-y_2=0$ and $I_3=x_3-y_3=0$, which is invariant under the dynamics of the joint system of Eqs.~(\ref{exA1}) and~(\ref{C1}).
Let us also define the complementary coordinates $J_i=y_i$, as discussed in Definition 3.

It is straightforward to check, according to definitions 2 and 3,  that ${\bf I}({\bf x},{\bf y})={\bf x}-{\bf y}$ is a three-dimensional non-singular semi-invariant that obeys the following ODEs system:
\begin{equation}
\frac{d{\bf I}}{dt}={\bf L}\left({\bf I},{\bf J}\right){\bf I},\;\;{\bf J}={\bf y},
\end{equation} 
where the matrices ${\bf L}$ and ${\bf L}_S$, according to definitions 5 and 6, are given by:
\begin{equation}
{\bf L}={\bf L}\left({\bf I},{\bf J}\right)= \left[
{\begin{array}{ccc}
 - \sigma  - c & \sigma  & 0 \\
r-J_3-I_3 & -1 &  -J_1-I_1 \\
J_2+I_2 & J_1+I_1 &  - b
\end{array}}
	\right],\hspace{5mm}
\label{matlorenz}
{\bf L}_S={\bf L}\left({\bf 0},{\bf J}\right)=
\left[
{\begin{array}{ccc}
 - \sigma  - c & \sigma  & 0 \\
r - {J}_{3} & -1 &  - {J}_{1} \\
{J}_{2} & {J}_{1} &  - b
\end{array}}
 \right].
\end{equation}

For a given initial condition ${\bf x}(t_0)={\bf x}_0$ and ${\bf y}(t_0)={\bf y}_0$ at initial time $t_0$, we define the impulsive vector field ${\bf G}_n$ and the respective times $t_n$ as follows - see Eq.~(\ref{formula16}):
\begin{eqnarray}
\label{eq49}&{\bf x}(t_n^+)={\bf y}(t_n)+\exp(A_n)\left({\bf x}(t_n^-)-{\bf y}(t_n)\right),
\;\;\displaystyle A_n=-\alpha\delta-\ln\left(\frac{||{\bf x}(t_n^-)-{\bf y}(t_n)||}{||{\bf x}(t_{n-1}^+)-{\bf y}(t_n)||}\right),&\nonumber\\
&t_{n+1}=t_n+\delta,\;\;t_1>t_0\;\;(n=1,2,\ldots,\infty),&
\end{eqnarray}
where $\alpha,\,\delta > 0$. 

Let us remark that $A_n$ in Eq.~(\ref{eq49}) corresponds to the exponent defined in Eq.~(\ref{formula20}), then from Eq.~(\ref{formula18}) and  Eq.~(\ref{formula21}) we have $B_n=-\alpha\delta$. Moreover, as $\Delta_n=t_n-t_{n-1}=\delta$, then the impulsive vector field defined in Eq.~(\ref{eq49}) fulfills the condition of proposition 4. Therefore, the trajectory of the impulsive system defined by equations (\ref{exA1}), (\ref{C1}) and (\ref{eq49}) converge to the synchronization surface $S$ for any initial condition. 

\begin{figure}[ht]
\begin{center}
	\includegraphics[scale=0.3]{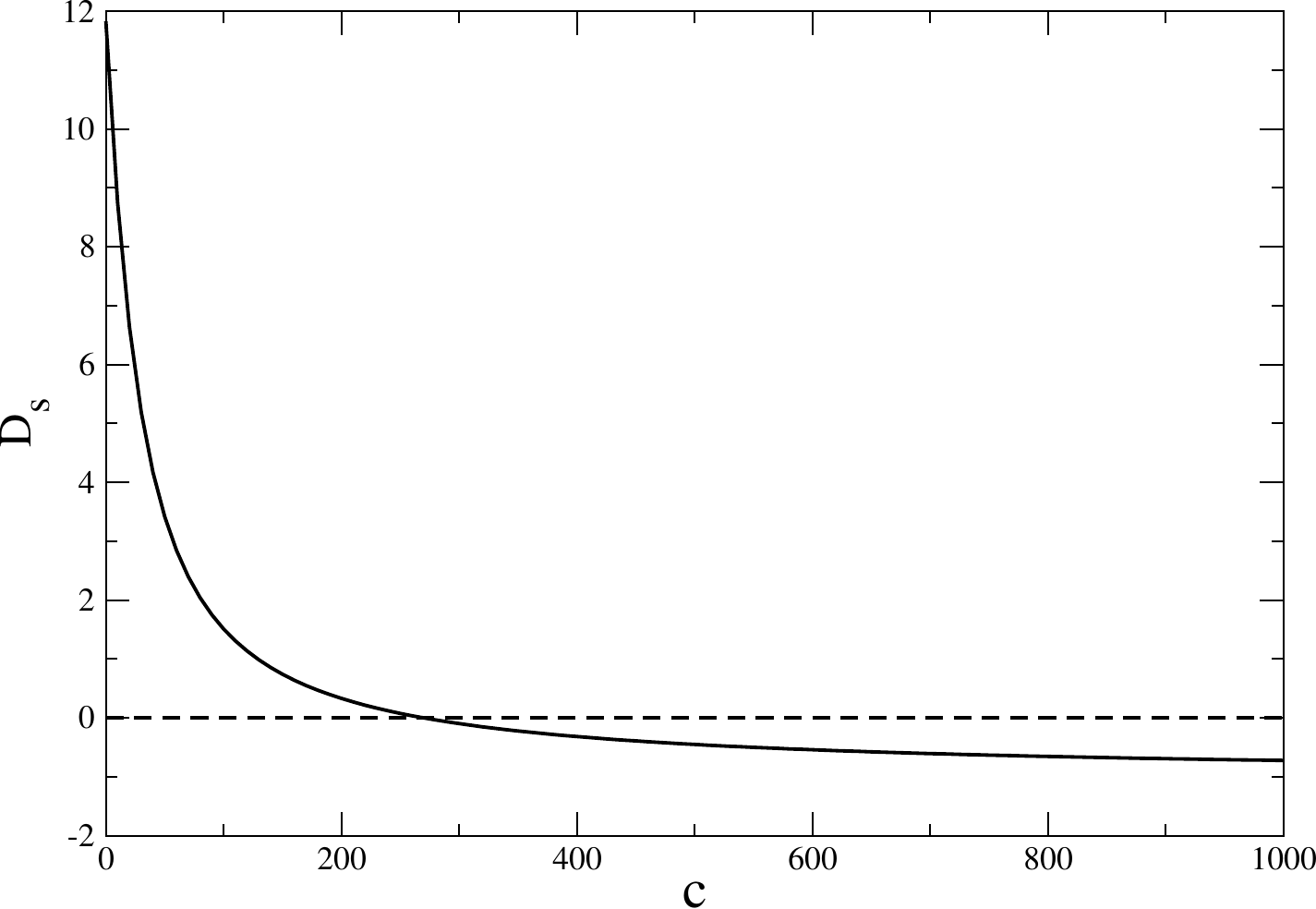}
\end{center}
\caption{Stability exponent $D_S$ as a function of the coupling parameter $c$.}
\label{fig2}
\end{figure}

Figure~\ref{fig2} shows the stability exponent $D_S$ as a function of the coupling
parameter $c$ in the absence of any forcing.
We see that there is a critical value $c_0>0$ for which the highest TLE of the synchronization surface
is negative for $c>c_0$.

To illustrate the impulsive control defined above  we present and application by choosing $b=8/3$, $\sigma=10$, $r=28$ and $c=5$ for the Lorenz coupled system. In this case, the stability exponent $D_S$ is positive. For the parameters that define the impulses, we choose $\delta=0.1$ and two values for $\alpha$:
$\alpha=0.1$ and $\alpha=0.4$.

\begin{figure}[ht]
\begin{center}
	\includegraphics[scale=0.3]{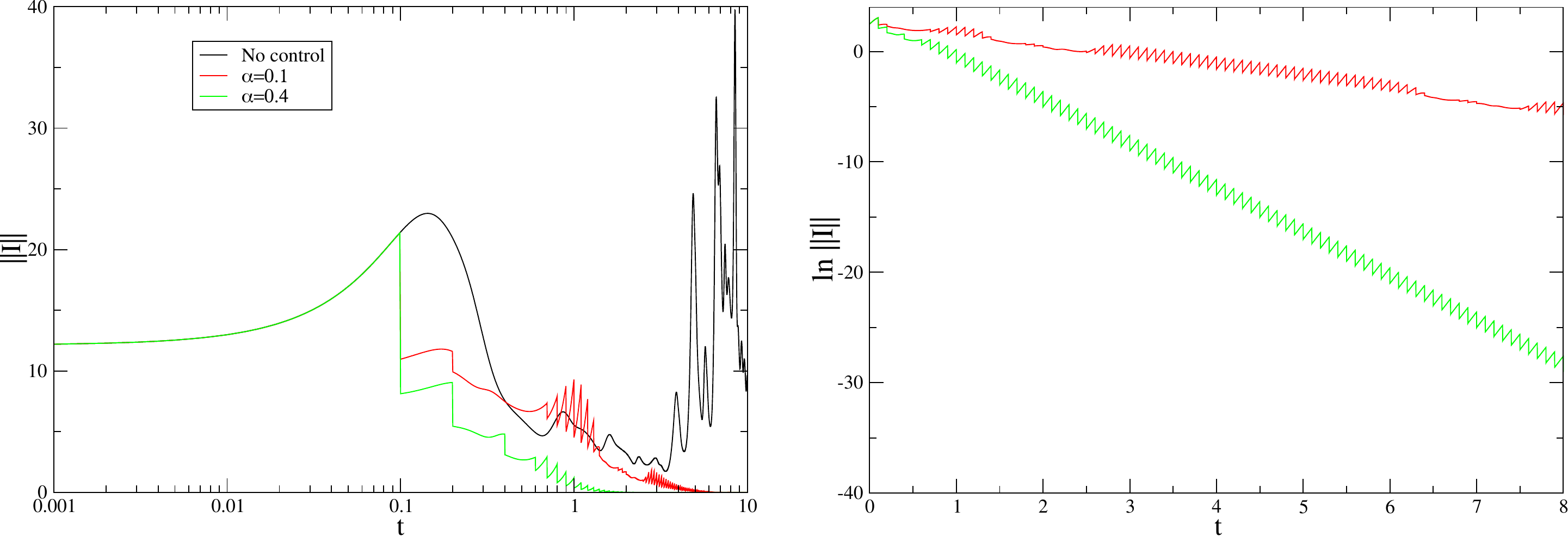}
\end{center}
\caption{Left panel: norm of $||{\bf I}||$ as a function of time for the cases without control
and for impulses with $\alpha=0.1$ and $\alpha=0.4$. Right panel:  $\log ||{\bf I}||$ for the cases with an impulsive control.}
\label{fig3}
\end{figure}

The left panel of Fig.~\ref{fig3} shows the time evolution of $||{\bf I}(t)||$ for the initial conditions $x_i=3$, $y_i=10$, $i=1,2,3$,
and the value of $\ln ||{\bf I}||$ is shown in the right panel of Fig.~\ref{fig3}.
Both control cases illustrated our results, putting in evidence the convergence of
the trajectories towards the synchronization manifold $S$, with a roughly exponential decay of $||{\bf I}(t)||$ with time.

\section{Parallel Impulsive Control}

According to Eq.~(\ref{formula17}), all impulsive control systems, illustrated in section \ref{sec5}, have impulsive vector field ${\bf G}_n$
such that ${\bf G}_n^J={\bf 0}$, which means that only the semi-invariant vector ${\bf I}$ are modified by the impulses, while the vector ${\bf J}$ is kept unchanged.
Indeed, in order to define an impulsive control system that satisfies the theorems 2 and 3, it is necessary that the semi-invariant vector ${\bf I}$ should be modified after each impulse.
However, it is interesting to think about impulsive vector fields that only changes the vector ${\bf J}$, which means to have $\Delta {\bf I}(t_n)={\bf G}_n^I={\bf 0}$ for all $n\in\mathbb{N}$. We call this type of control Parallel Impulsive Control.

The main idea, based on the integration of Eq.~(\ref{formula10}) in proposition 3, is to define a parallel control that satisfies the following condition:
\begin{equation}
\label{paralelo1}
\lim_{t\rightarrow\infty}\int_{t_0}^t \left<{\bf i},{\bf H}{\bf i}\right>_{t'}dt'=-\infty 
\;\;\Rightarrow\;\;\lim_{t\rightarrow\infty}||{\bf I}(t)||=||{\bf I}(t_0)||\exp\left(\lim_{t\rightarrow\infty}\int_{t_0}^t \left<{\bf i},{\bf H}{\bf i}\right>_{t'}dt'\right)=0. 
\end{equation}    
Let us observe that for a parallel control, the scalar product, which is integrated in Eq.~(\ref{paralelo1}), is modified by impulses only by changes in the matrix ${\bf H}$, since the versor ${\bf i}={\bf I}/{||{\bf I}||}$ is kept unchanged.

Le us apply our approach of parallel impulsive control to a
SEIR epidemiological model for measles~\cite{viralinfec}, with the addition
of a vaccinated population variable $V$~\cite{ref52,ref53}, where
the possibility of optimizing a vaccination protocol is particularly relevant in real world situations with finite resources. 
In the absence of vaccination the equations governing the model are:
\begin{eqnarray}\nonumber
	 & & \frac{dV}{dt} = 0, \hspace{5mm}\frac{dS}{dt} = -\rho\frac{S I}{N},\nonumber\\
         & & \frac{dE}{dt} = \rho\frac{S I}{N} - \sigma E,\nonumber\\
	 & & \frac{dI}{dt} = \sigma E - \gamma I,\hspace{5mm}\frac{dR}{dt} = \gamma I,
\label{seir}
\end{eqnarray}
where $V$, $S$, $E$, $I$ and $R$ are the proportions with respect to the total population $N$ of
vaccinated, susceptible, exposed, infected and recovered
individuals, respectively, $N=(V+S+E+I+R)=1$ and $\rho$ the transmission rate, $\sigma$ the incubation rate and $\gamma$ the recovery rate.

The semi-invariant vector (see definitions 2 and 3) for the system in Eq.~(\ref{seir}) is given by ${\bf I}=(E,I)^t$, the matrices ${\bf L}$ and ${\bf H}$ (see definition 5) are given respectively by
\begin{equation}
{\bf L}=\left[\begin{array}{cc} -\sigma &  \rho S \\ \sigma & -\gamma \end{array}\right],\;\;
{\bf H}=\left[\begin{array}{cc} -\sigma &  \displaystyle \frac{\sigma+\rho}{2} S \\ \displaystyle\frac{\sigma+\rho S}{2} & -\gamma \end{array}\right],
\end{equation}
and, from Eq.~(\ref{formula6}), the vector ${\bf J}$  is given by ${\bf J}=(S,V,R)^t$. Let us observe that the matrix ${\bf H}$ depends only on the variable $S$.

In what follows the variable $V$ is used as an intervention variable such that the time evolution converges to the disease free
invariant surface defined as
\begin{equation}
        \{(E,I,S,V,R)\;|\;E=0,I=0\}.
        \label{Sarbo}
\end{equation}
We then considered the parallel impulsive system defined as
\begin{equation}
\label{paralelo2}
\begin{array}{c}
\begin{array}{c}
	S(t_{n}^{+}) = \exp(B_n)S(t_{n}^{-}), \\ \\
	V(t_{n}^{+}) = V(t_{n}^{-}) + [1-\exp(B_n)]S(t_{n}^{-}),
\end{array}
\;\;t_{n+1}=t_n+\kappa_{\mbox{\small eff}}(t_n-t_0),\;\; t_1>t_0\;\;(n=1,2,\ldots,\infty),\;\; \kappa_{\mbox{\small eff}}>0,\\ \\
\displaystyle
B_{n} = -\ln\frac{||{\bf I} (t^{-}_{n})||}{||{\bf I} (t^{+}_{n-1})||}-\alpha\left(t_{n}-t_{n-1}\right)\;\;\;
	\textrm{if}\;\;\;||{\bf I}(t_{n+1})||>||{\bf I}(t_{n+1})||\;\;\left(||{\bf I}||=\sqrt{I^2 + E^2}\right),
\end{array}
\end{equation}
where at each impulse a proportion of the susceptible population is vaccinated (moving from the susceptible to the vaccinated class $V$):

We consider here the initial conditions $V(0)=0.3$, $E(0)=0$, $I(0)=2\times 10^{-4}$, and parameter values $\gamma=365/7$ (average recovery time of 7 days), $\sigma=365/8.5$ (average incubation time of 8.5 days) and $\rho=114.715$~\cite{ref52,ref53}. The condition given in Eq.~(\ref{paralelo1}) is satisfied for any sequence of impulses, because the SEIR system in Eq.~(\ref{seir}), for any initial value of $S$, converges to a some disease free invariant surface as defined in Eq.~(\ref{Sarbo}). The main objective here is not to constraint the convergence toward the invariant surface, but to control the infection spread velocity with a minimal number of impulses. 

\begin{figure}[ht]
\centering
\includegraphics[scale=0.3]{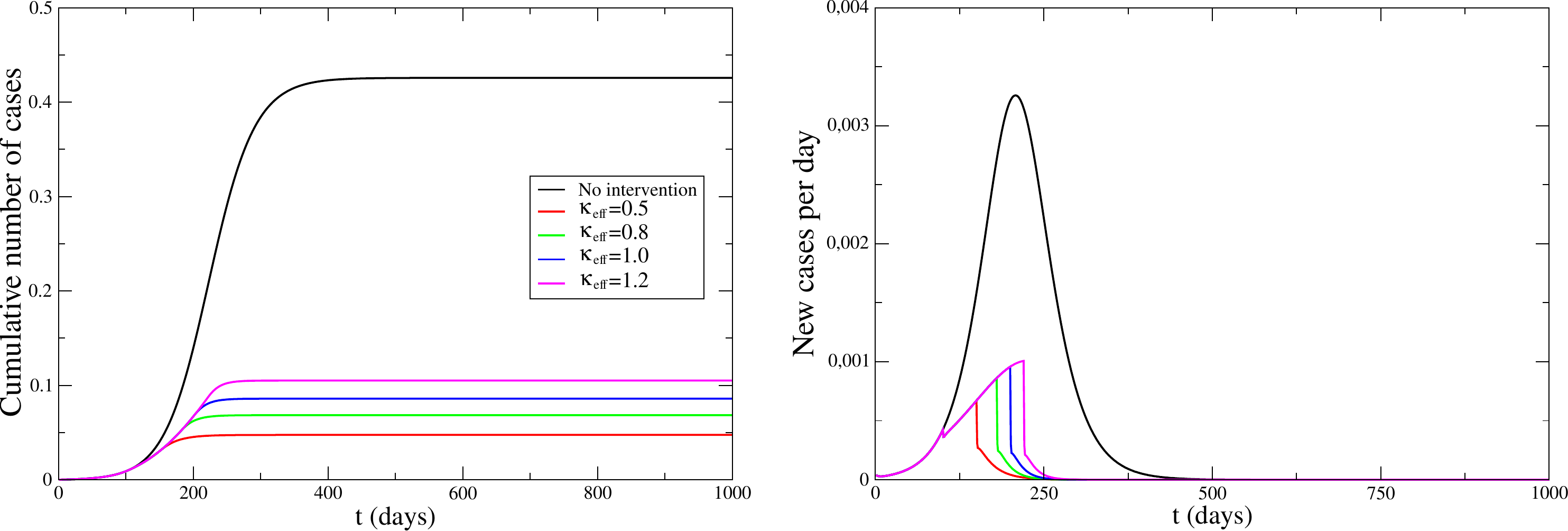}
\caption{Left panel: evolution of the cumulative number of cases as a proportion of the total population.
	Right panel: epidemic curves (new cases per day) as a proportion to the total population.}
\label{sarampo_original}
\end{figure}

The first impulse is chosen to occur at $t_1=$ 100 (days) ($t_0=0$) and the convergence speed control parameter
as $\alpha=0.002$.
Figure- \ref{sarampo_original} shows the cumulative number of cases for the  outbreak of measles for $\kappa_{\mbox{\small eff}}=0.5, 0.8, 1.0, 1.2$
for a period of 3 years, alongside the natural evolution in the absence of intervention.
The number of measles cases rapidly converges to the disease free invariant surface for all the impulsive vaccination considered,
with a small number of interventions. This is particularly important for vaccination designs with a limited amount of resources.
The present approach can be easily extended to other epidemiological models describing diseases with an existing vaccine.

\section{Concluding Remarks}
\label{sec6}

We presented here a novel approach of impulsive control for the convergence of any trajectory towards an invariant manifold
of a system of ordinary differential equations, and stated the suitable conditions for its realization.
For different choices of impulsive interventions, as defined by the choice of the impulsive vector field,
the stability condition of the control is expressed by Theorem 3. This results leads to
a great flexibility for choosing control parameters that are convenient in practical situations.
Our approach is feedback-adaptive as it requires the determination of the values of ${\bf I}$
at the successive intervention times $t_n^{-}$ and $t_{n-1}^{+}$.
A convergence speed parameter was defined to control convergence speed.

The knowledge of the stability exponent is not necessary in our approach, provided the intervals
between two successive impulses are upper bounded.
In our applications for the Lorenz system 
the time intervals between two successive impulses grow with time, illustrating the possibility to carry out 
a control with interventions becoming less frequent. This is particularly relevant in practical situations where the cost of the intervention is important.
Finally, we also observe that our approach is such that any divergence caused by an invariant object belonging to the synchronization surface
is canceled out by the impulses. This avoid the desynchronization observed for some continuous control approaches~\cite{r5}.

Finally, we observe that theorem 2 establishes, for an impulsive system having a semi-invariant with well defined stability exponent, a necessary and sufficient condition that assure the stability for the invariant surface defined by the semi-invariant, while theorem 3 establishes only a sufficient condition for a semi-invariant without the assumption about the stability exponent existence. However, both theorems can not be applied for parallel impulsive systems and, in this case, the general condition to guarantee the stability is given by Eq.~(\ref{paralelo1}).  A more detailed study on necessary and sufficient conditions for parallel impulsive systems, establishing detailed condition for the interval between two impulses, must be developed and, we hope, it  will be the subject of  a future work.


\end{document}